\documentstyle[12pt, leqno]{article}
\newtheorem{theorem}{Theorem}
\newtheorem{definition}[theorem]{Definition}
\newtheorem{proposition}[theorem]{Proposition}
\newtheorem{corollary}[theorem]{Corollary}
\newtheorem{lemma}[theorem]{Lemma}

\def \proof {\noindent {\bf Proof.}\ \ }
\def \endproof {{\mbox{}\nolinebreak\hfill\rule{2mm}{2mm}\par\medbreak} }

\def \ll {\langle}
\def \rr {\rangle}
\def \N {{\bf N}}
\def \R {{\bf R}}
\def \a {\alpha}
\def \e {\varepsilon}
\def \d {\delta}
\def \l {\lambda}
\def \s {\sigma}
\def \t {\tau}

\def \R {{\bf R}}
\def \ra {\rightarrow}

\def \dist {{\rm dist}}
\def \span {{\rm span}}
\def \rank {{\rm rank}}
\def \im {{\rm Im}}
\def \lo {\mbox{\large $0$}}
\def \halfn {{\lceil n/2 \rceil}}
\def \nm {{\sqrt \frac{n}{m}}}
\def \mn {{\sqrt \frac{m}{n}}}

\begin{document}
\title {Subsequences of frames}
\author {R. Vershynin\footnote{Department of Mathematics, University of Missouri, 
                               Columbia, MO 65211, USA (vershynin@yahoo.com)} }
\date{February 13, 1999}
\maketitle

\begin{abstract} 
  Every frame in Hilbert space contains a subsequence equivalent to an orthogonal basis.
  If a frame is $n$-dimensional then this subsequence has length $(1 - \e) n$.
  On the other hand, there is a frame which does not contain bases with brackets.
\end{abstract}


\section{Introduction}                             

The notion of frame goes back to R.Duffin and A.Schaeffer \cite{D-S} and was studied 
extensively since then with relation to nonharmonic Fourier analysis, see \cite{He}.
From a geometrical point of view, a frame in a Hilbert space $H$ is the image
of an orthonormal basis in a larger Hilbert space under an orthogonal projection
onto $H$, up to equivalence \cite{Ho}
(the equivalence constant is called the frame constant).
Since frames have nice representation properties (see \cite{D-S}, \cite{A}), 
much attention was paid to their subsequences that inherit these properties.
The most interesting questions arise about subsequences equivalent to an orthogonal
basis \cite{Ho}, \cite{S}, \cite{C1}, \cite{C-C1}.
P.Casazza \cite {C2} proved that, given an $\e > 0$, 
any $n$-dimensional frame whose norms are well bounded below
contains a subsequence of length $(1 - \e) n$ equivalent to an orthogonal basis
(the constant of equivalence does not depend of $n$). 

In the present paper this is proved for all frames, 
without restrictions on norms of the elements. 
If a frame is $n$-dimensional then it contains a subsequence 
of length $(1 - \e) n$ which is $C$-equivalent to an orthogonal basis. 
Here $C$ depends only on the frame constant and $\e$.
To put the result in other words, orthogonal projections in Hilbert space
preserve orthogonal structure in almost whole range. 
Namely, any orthogonal projection $H$ of an orthogonal basis 
contains a subset of cardinality $(1 - \e) \rank (P)$ 
which is $C(\e)$-equivalent to an orthogonal system. 
This is proved in Section \ref{secfinitedim}.

An infinite dimensional version of this result is considered in Section \ref{basicsubs}.
Every infinite dimensional frame has an infinite subsequence equivalent to an 
orthogonal basis.
However, for some frames this subsequence can not be complete, 
as was shown by K.Seip \cite{S} and P.Casazza and O.Christensen \cite{C-C2}.
This result is generalized in Section \ref{secwithoutbrackets} 
by constructing a frame which does not contain bases with brackets. 
So our frame $(x_j)$ is "asymptotically indecomposable" in the following sense.
If $(y_j)$ is any complete subsequence of $(x_j)$, then the distance 
from $\span(y_j)_{j \le n}$ to $\span(y_j)_{j > n}$ tends to zero as $n \ra \infty$.

In the rest of this section we recall standard definitions and simple
known facts about frames. 
In what follows, $H$ will denote a separable Hilbert space,
finite or infinite dimensional. 
Absolute constants will be denoted by $c_1, c_2, \ldots$.
A sequence $(x_j)$ in $H$ is called a {\em frame} 
if there exist positive numbers $A$ and $B$ such that
$$
A \|x\|^2  \le  \sum_j | \ll x, x_j \rr |^2  \le  B \|x\|^2
\ \ \ \ \mbox{for $x \in H$.}
$$
The number $(B/A)^{1/2}$ is called a {\em constant} of the frame.
We call $(x_j)$ a {\em tight frame} if $A = B = 1$.

Two sequences $(x_j)$ and $(y_j)$ in possibly different Banach spaces 
are called {\em equivalent} if there is an isomorphism 
$T : [x_j] \ra [y_j]$ such that $Tx_j = y_j$ for all $j$. 
Here $[x_j]$ denotes the closed linear span of $(x_j)$.
Let $c  =  \|T\| \|T^{-1}\|$ then the sequences $(x_j)$ and $(y_j)$ 
are called {\em $c$-equivalent}.

\medskip

The next observation (see \cite{Ho}) allows to look at frames 
as at projections of the canonical vector basis $(e_j)$ in $l_2$.

\begin{proposition}                                                   
        \label{view}
  Let $(x_n)_{n=1}^m$ be a frame in $H$ with constant $c$, 
  where $m$ can be equal to infinity.
  Then there is an orthogonal projection $P$ in $l_2^m$ such that  
  $(x_n)$ is $c$-equivalent to $(P e_n)$.
  Conversely, if $P$ is an orthogonal projection in $l_2^m$ onto a
subspace $H$, 
  then $(P e_n)_{n=1}^m$ is a tight frame in $H$.
\end{proposition}

%

\begin{corollary}                                 \label{frameistight}
  Let $(x_n)$ be a frame with constant $c$. 
  Then $(x_n)$ is $c$-equivalent to a tight frame. 
\end{corollary}

Now we present another view at frames. 
We can regard them as the columns of a row-orthogonal matrix (either
finite or infinite).

\begin{lemma}                                      \label{columnrow}
  Let $n, m \in \N \cup \infty$ and $A$ be an $n \times m$ matrix
  whose rows are orthonormal. 
  Then the columns of $A$ form a tight frame in $l_2^n$.
  
  Conversely, let $(x_j)_{j=1}^m$ be a frame in $H$. 
  Then there exists an $n \times m$ matrix $A$ with $n = \dim H$
  whose rows are orthonormal and such that the columns 
  form a tight frame equivalent to $(x_j)$.
\end{lemma}

\proof
If $A$ is as above then $A^*$ acts as an isometric embedding 
of $l_2^n$ into $l_2^m$.
Then $A$ acts as a quotient map in a Hilbert space, and we can regard
it 
as an orthogonal projection. On the other side, the columns of $A$
are equal to $Ae_j$. Proposition \ref{view} finishes the proof of the
first statement.
The converse can also be proved by this argument.
\endproof

\begin{lemma}                                                   \label{square}
  Let $(x_j)$ be a tight frame in $H$. 
  Then $\sum_j \|x_j\|^2  = \dim H$ 
  (which is possibly equal to infinity).
\end{lemma}

\proof 
By Proposition \ref{view} we may assume that $H$ is a subspace of $l_2$ 
and $x_j = P e_j$, where $P$ is the orthogonal projection in $l_2$
onto $H$.
Then the Hilbert-Schmidt norm 
$\|P\|_{\rm HS}  =  ( \sum_j \|x_j\|^2 )^{1/2}$. 
On the other hand, $\|P\|_{\rm HS}  =  (\dim H)^{1/2}$.
\endproof

%
%
%

\section{Finite dimensional frames}                              \label{secfinitedim}

%
%
%
%
%
%

In this section we prove 

\begin{theorem}                                          \label{quantitative}
There is a function $h : \R_+ \ra \R_+$ such that the following holds. 
Suppose  $(x_j)$ is an $n$-dimensional frame with constant $c$. 
Then for every $\e > 0$ there is a set of indices $\s$ 
with $|\s| > (1-\e) n$ such that 
the system $(x_j)_{j \in \s}$ 
is $C$-equivalent to an orthogonal basis, 
where $C = h(\e) c$. 
\end{theorem}

We will need a result of A.Lunin on norms of restriction of operators
onto coordinate subspaces \cite{L} (for improvements see \cite{K-Tz}).

\begin{theorem} (A.Lunin).                                 \label{Lunin}
  Let $T : l_2^m  \ra  l_2^n$ be a linear operator.
  Then there is a set $\s  \subset  \{1, \ldots, m\}$ with $|\s| = n$
  such that 
  $$
  \| T |_{\R^\s} \|  \le  c_1 \nm \|T\|.
  $$
\end{theorem}

Given an $h>0$, a system of vectors $(x_j)$ in a Hilbert space is
called 
{\bf $h$-Hilbertian} if
$$
\Big\| \sum_j a_j x_j \Big\|   \le   h \Big( \sum_j |a_j|^2 \Big)^{1/2}
$$
for all sequences of scalars $(a_j)$.
Then Theorem \ref{Lunin} can be reformulated as follows.
Suppose $(x_j)_{1 \le j \le m}$ is a $1$-Hilbertian system 
in $l_2^n$. Then there is a set $\s  \subset  \{1, \ldots, m\}$ 
with $|\s| = n$ such that 
$(\mn x_j)_{j \in \s}$ is $c_1$-Hilbertian.

Next, we will use a result of J.Bourgain and L.Tzafriri 
on invertibility of large submatrices \cite{B-Tz} Theorem 1.2:

\begin{theorem} (J.Bourgain, L.Tzafriri).                             
      \label{BourgainTzafriri}  
  Let $T : l_2^n  \ra  l_2^n$ be a linear operator 
  such that $\|T e_j\| = 1$ for all $j$. 
  Then there is a set $\s  \subset  \{1, \ldots, n\}$ 
  with $|\s|  \ge  c_2 n / \|T\|^2$ such that
  $$
  \|Tx\|  \ge  c_2 \|x\|  \ \ \ \ \mbox{for every $x \in \R^\s$}.
  $$
\end{theorem}

Given a $b>0$, a system of vectors $(x_j)$ in a Hilbert space is called 
{\bf $b$-Besselian} if
$$
b \Big\| \sum_j a_j x_j \Big\|   \ge   \Big( \sum_j |a_j|^2 \Big)^{1/2}
$$
for all sequences of scalars $(a_j)$.
Then Theorem \ref{BourgainTzafriri} can be reformulated as follows. 
Suppose $(x_j)_{1 \le j \le n}$ is an $h$-Hilbertian system 
in $l_2^n$ and $\|x_j\| \ge \a$ for all $1 \le j \le n$.
Then there is a set $\s  \subset  \{1, \ldots, n\}$ 
with $|\s|  \ge  c_2 (\a / h)^2 n$ such that 
the system $(\a^{-1} x_j)_{j \in \s}$ is $c_3$-Besselian.

Clearly, every tight frame is $1$-Hilbertian.

\begin{lemma}                                                       
\label{sizeoftau}
  Let $(y_j)_{1 \le j \le m}$ be a tight frame in $l_2^n$
  with $\|y_j\| = \nm$ for all $j$.
  Let $P$ be a $k$-dimensional orthogonal projection in $l_2^n$.
  Then for $\d > 0$
  $$
  \Big| \Big\{  j : \| (I - P) y_j \|  \ge  \d \nm \Big\} \Big|
  \ge  \Big( 1 - \d^2 - \frac{k}{n} \Big) m.
  $$
\end{lemma}

\proof   
Let $\t = \Big\{  j : \| (I - P) y_j \|  \ge  \d \nm \Big\}$.
Since $((I - P)y_j)_{1 \le j \le m}$ is a tight frame in an
$(n-k)$-dimensional space $(I - P) l_2^n$, Lemma \ref{square} yields
\begin{eqnarray*}
n - k  &=&   \sum_{j=1}^m \|(I - P) y_j\|^2
   \le   \sum_{j \in \t} \|y_j\|^2  +  \sum_{j \in \t^c} \|(I -
P)y_j\|^2      \\
  &\le&   |\t| \cdot (n/m)  +  m \cdot \d^2 (n/m)
   =    (|\t| / m + \d^2) n.
\end{eqnarray*}
The required estimate follows.
\endproof

\vspace{0.5cm}
Now we proceed to the proof of Theorem \ref{quantitative}.
As in P.Casazza's proof \cite{C2}, the set $\s$ will be constructed by 
an iteration procedure. Our proof consists of several parts.

{\bf I. Splitting.}  \ \
By Corollary \ref{frameistight}, we may assume that 
the frame $(x_j) \subset l_2^n$ is tight and all of its terms are nonzero. 
First we will split $(x_j)$ to get almost equal norms of the terms. 
Note that if we substitute any member $x_j$ of the frame by $k$ elements
$x_j / \sqrt{k}, \ldots, x_j / \sqrt{k}$, we will still get a tight
frame. 
Fix a $\nu > 0$.
Splitting each element $x_j$ as above, we can obtain a new tight frame
$(y_j)_{1 \le j \le m}$ such that

(i) elements of $(y_j)$ are multiples of the ones from $(x_j)$;

(ii) there is a $\l > 0$ such that 
     $\l  \le  \|y_j\|  \le  (1 + \nu) \l$ for all $j = 1, \ldots, m$.

\noindent The constant $\l$ be evaluated using Lemma \ref{square}:
$$
(1 + \nu)^{-1} \nm  \le  \|y_j\|  \le  (1 + \nu) \nm
   \ \ \ \ \mbox{for $j = 1, \ldots, m$}.
$$
Clearly, it is enough to prove the theorem for $(y_j)$ instead of
$(x_j)$.
We can choose the parameter $\nu =\nu(\e) > 0$ arbitrarily small. 
To make the proof more readable, we simply assume that $\nu = 0$
which is a slight abuse of rules. The reader will easily adjust the
arguments
to the general case. So we have
$$
\|y_j\|  =  \nm,      \ \ \ \ j=1, \ldots, m.
$$
We can also assume that $(\e / 2) m  \ge  n$.

{\bf II. Iterative construction.} \ \ 
Let $\d = \sqrt{\e / 2}$.

{\em Step 1.} \ 
Set $\t_0 = \{1, \ldots, m\}$.
The system $(y_j)_{j \in \t_0}$ is $1$-Hilbertian. 
Lunin's theorem yields the existence of a set $\s'_1 \subset \t_0$ 
with $|\s'_1| = n$ such that
$$
\mbox{the system $(\mn y_j)_{j \in \s'_1}$ is $c_1$-Hilbertian.}
$$ 
Note that $\| \mn y_j \| = 1$ for $j \in \s'_1$. 
Then Bourgain-Tzafriri's theorem gives us a set 
$\s_1 \subset \s'_1$ with $|\s_1|  \ge  (c_2 / c_1^2) n$ such that
$$
\mbox{the system $(\mn y_j)_{j \in \s_1}$ is $c_3$-Besselian.}
$$ 
So we have already found a subsequence $(y_j)_{j \in \s_1}$ 
of length proportional to $n$ which is well equivalent
to an orthogonal basis. If $|\s_1|  \ge  (1 - \e) n$, 
then we are done and stop here. Otherwise proceed to the next step. 

{\em Step 2.} \  
Let $P_1$ be the orthogonal projection in $l_2^n$ 
onto $[y_j]_{j \in \s_1}$.
Let 
$$
\t_1  =  \Big\{ j  : \| (I - P_1) y_j \|  \ge  \d \nm \Big\}.
$$
Clearly, $\t_1  \subset  \s_1^c$.
By Lemma \ref{sizeoftau}
$$                               
  |\t_1|  \ge  \Big( 1 - \d^2 - \frac{|\s_1|}{n} \Big) m.
$$
As $|\s_1|  <  (1 - \e) n$, 
$$
|\t_1|  >  \Big( 1 - \d^2 - (1 - \e) \Big) m
        =  (\e / 2) m.
$$
The system $(y_j)_{j \in \t_1}$ is $1$-Hilbertian 
and $|\t_1|  \ge  n$ by the choise of $m$.
Lunin's theorem yields the existence of a set $\s'_2 \subset \t_1$ 
with $|\s'_2| = n$ such that
$$
\mbox{the system 
  $( \sqrt{\frac{|\t_1|}{n}} y_j )_{j \in \s'_2}$ 
  is $c_1$-Hilbertian.}
$$ 
Then the system $(\sqrt{\frac{|\t_1|}{n}} (I - P_1) y_j)_{j \in \s'_2}$ 
is also $c_1$-Hilbertian.
By the definition of $\t_1$, it has not too small norms:
$$
\Big\| \sqrt{\frac{|\t_1|}{n}} (I - P_1) y_j \Big\|
\ge  \d \sqrt{\frac{|\t_1|}{m}},  \ \ \ \ j \in \s'_2.
$$
Then Bourgain-Tzafriri's theorem gives us a set $\s_2 \subset
\s'_2$ 
with 
$$
|\s_2|  \ge  c_2 \Big( \d^2 \frac{|\t_1|}{m} / c_1^2 \Big) n
\ge  (c_2 / c_1^2) \d^2 \Big( (1 - \d^2) n - |\s_1| \Big)
$$
such that
$$
\mbox{the system 
  $( \mn (I - P_1) y_j )_{j \in \s_2}$
  is $(c_3 \d^{-1})$-Besselian.}
$$
If $|\s_1| + |\s_2|  \ge  (1 - \e) n$, then we stop here. 
Otherwise proceed to the next step.
 
{\em Step $k+1$.} \ We assume that the sets $\s_1, \ldots, \s_k$ 
are already constructed and
\begin{equation}                                                  \label{notyet}
  \sum_{i=1}^k |\s_i|  <  (1 - \e) n.
\end{equation}
Let $P_k$ be the orthogonal projection in $l_2^n$ 
onto $[y_j]_{j \in \s_1 \cup \ldots \cup \s_k}$.
Let 
$$
\t_k  =  \Big\{ j  : \| (I - P_k) y_j \|  \ge  \d \nm \Big\}.
$$
Clearly, $\t_k  \subset  (\s_1 \cup \ldots \cup \s_k)^c$.
By Lemma \ref{sizeoftau}
$$
|\t_k|  \ge  \Big( 1 - \d^2 - \frac{ \sum_{i=1}^k |\s_i| }{n} \Big) m.
$$
By (\ref{notyet})
$$
|\t_k|  >  \Big( 1 - \d^2 - (1 - \e) \Big) m
        =  (\e / 2) m.
$$
The system $(y_j)_{j \in \t_k}$ is $1$-Hilbertian 
and $|\t_k|  \ge  n$ by the choise of $m$.
Lunin's theorem yields the existence of a set $\s'_{k+1} \subset \t_k$ 
with $|\s'_{k+1}| = n$ such that
$$
\mbox{the system 
  $( \sqrt{\frac{|\t_k|}{n}} y_j )_{j \in \s'_{k+1}}$ 
  is $c_1$-Hilbertian.}
$$ 
Then the system $(\sqrt{\frac{|\t_k|}{n}} (I - P_k) y_j)_{j \in \s'_{k+1}}$ 
is also $c_1$-Hilbertian.
By the definition of $\t_k$, it has not too small norms:
$$
\Big\| \sqrt{\frac{|\t_k|}{n}} (I - P_k) y_j \Big\|
\ge  \d \sqrt{\frac{|\t_k|}{m}},  \ \ \ \ j \in \s'_{k+1}.
$$
Then Bourgain-Tzafriri's theorem gives us a set $\s_{k+1} \subset \s'_{k+1}$ 
with 
\begin{equation}                                             
\label{sizeofsigma}
  |\s_{k+1}|  \ge  c_2\Big( \d^2 \frac{|\t_k|}{m} / c_1^2 \Big) n
  \ge  (c_2 / c_1^2) \d^2 \Big( (1 - \d^2) n - \sum_{i=1}^k|\s_i| \Big)
\end{equation}
such that
$$
\mbox{the system 
  $( \mn (I - P_k) y_j )_{j \in \s_{k+1}}$
  is $(c_3 \d^{-1})$-Besselian.}
$$
If $\sum_{i=1}^{k+1} |\s_i|  \ge  (1 - \e) n$, then we stop here. 
Otherwise proceed to the next step.

{\bf III. When we stop.} \ \ 
Let $k_0$ be the number of the last step, that is the smallest integer
such that
$$
\sum_{i=1}^{k_0} |\s_i|  \ge  (1 - \e) n.
$$
We claim that such $k_0$ exists 
and there is a function $K(\e)$ such that $k_0  \le  K(\e)$.
Indeed, let $K(\e) = [4 c_1^2 c_2^{-1} \e^{-2}] + 2$. 
If the claim were not true, then 
$$
\sum_{i=1}^k |\s_i|  <  (1 - \e) n     \ \ \ \ \mbox{for $k = 1, \ldots, K(\e)$}.
$$
Then by (\ref{sizeofsigma}) for all $k = 2, \ldots, K(\e)$
\begin{eqnarray*}
|\s_k|  &\ge&  (c_2 / c_1^2) \d^2 \Big( (1 - \d^2) - (1 - \e) \Big) n               \\
        & = &  (c_2 / c_1^2) (\e^2 / 4) n.
\end{eqnarray*}
Thus 
$$
\sum_{i=1}^{K(\e)} |\s_i|  
\ge  (K(\e) - 1) \cdot (c_2 / c_1^2) (\e^2 / 4) n  \ge  n.
$$
This contradiction proves the claim.

Now set $\s = \s_1 \cup \ldots \cup \s_{k_0}$, 
then $|\s| > (1 - \e) n$.
To complete the proof of the theorem, it remains to check 
that the system $(\mn y_j)_{j \in \s}$ is well equivalent 
to an orthonormal basis.

{\bf IV. Equivalence to the orthogonal basis within blocks $\s_k$.} \ \ 
Recall that for every $k < k_0$ the size of $\t_k$ is comparable with $m$,
namely $|\t_k|  \ge  (\e / 2) m$.
Then we conclude from the construction the existence of functions 
$c_1(\e)$ and $c_2(\e)$ such that for every $k = 1, \ldots, k_0$
\begin{equation}                                                  
\label{hilbertian}
  \mbox{the system 
    $(\mn y_j )_{j \in \s_k}$ 
    is $c_1(\e)$-Hilbertian,}
\end{equation}
\begin{equation}                                                  
\label{besselian}
  \mbox{the system 
  $(\mn (I - P_{k-1}) y_j )_{j \in \s_k}$
  is $c_2(\e)$-Besselian.}
\end{equation}

{\bf V. The system $(\mn y_j)_{j \in \s}$ is $h$-Hilbertian 
for some function $h = h(\e)$.} \ \ 
Indeed, fix scalars $(a_j)_{j \in \s}$ such that $\sum_{j \in \s}
|a_i|^2 = 1$.
Then 
\begin{eqnarray*}
\Big\| \sum_{j \in \s} a_j \Big( \mn y_j \Big) \|   
  &\le&  \sum_{k=1}^{k_0} \Big\| \sum_{j \in \s_k} a_j \Big( \mn y_j
\Big) \|    \\
  &\le&  \sqrt{k_0} \left( 
      \sum_{k=1}^{k_0} \Big\| \sum_{j \in \s_k} a_j \Big( \mn y_j
\Big) \Big\|^2
      \right)^{1/2}    \\ 
  &\le&  \sqrt{k_0} \; c_1(\e) \left(
      \sum_{k=1}^{k_0} \sum_{j \in \s_k} |a_j|^2 
      \right)^{1/2}         \ \ \ \ \mbox{by (\ref{hilbertian})}   \\ 
  & = &   \sqrt{K(\e)} \; c_1(\e).
\end{eqnarray*}

{\bf VI. The system $(\mn y_j)_{j \in \s}$ is $b$-Besselian 
for some function $b = b(\e)$.} \ \ 
We follow P.Casazza \cite{C2}. Choose $r = r(\e) > 2$ large enough 
(to be specified later). Let $a = a(\e) > 0$ be such that 
$r^{k_0+1} a  <  1$.
Fix scalars $(a_j)_{j \in \s}$ such that $\sum_{j \in \s} |a_j|^2 = 1$.
Suppose
\begin{equation}                                                   
\label{iolargest}
\mbox{$1 \le k' \le k_0$ is the largest so that} \ \ 
  \Big( \sum_{j \in \s_{k'}} |a_j|^2 \Big)^{1/2}
  \ge  r^{k_0 - k'} a.
\end{equation}
Such $k'$ must exist, otherwise
\begin{eqnarray*}
\Big( \sum_{j \in \s} |a_j|^2 \Big)^{1/2}
  &\le&   \sum_{k=1}^{k_0} \Big( \sum_{j \in \s_k} |a_j|^2 \Big)^{1/2}
\\
  &\le&   \sum_{k=1}^{k_0} r^k a
    \le  r^{k_0 + 1} a    <   1,
\end{eqnarray*}
contradicting the choise of $a$.
We have
\begin{eqnarray*}
\Big\| \sum_{j \in \s} a_j \Big( \mn y_j \Big) \Big\|  
  &\ge&    \Big\| \sum_{k=1}^{k'} \sum_{j \in \s_k} a_j \Big( \mn y_j \Big) \Big\|         
         - \sum_{k=k' + 1}^{k_0} \Big\|\sum_{j \in \s_k} a_j \Big( \mn y_j \Big) \Big\|        \\
  &\ge&    \Big\| (I - P_{k'-1}) \sum_{k=1}^{k'} \sum_{j \in \s_k} a_j \Big( \mn y_j \Big) \Big\| -     \\
  & & \ \ \ \ \ \ \  - c_1(\e) \sum_{k=k' + 1}^{k_0} \Big( \sum_{j \in \s_k} |a_j|^2 \Big)^{1/2}
         \ \ \ \ \mbox{by (\ref{hilbertian})}                      \\
  &\ge&    \Big\|  \sum_{j \in \s_{k'}} a_j \Big( \mn (I - P_{k' - 1}) y_j \Big) \Big\|         
         - c_1(\e) \sum_{k=k' + 1}^{k_0} r^{k_0 - k} a
         \ \ \ \ \mbox{by (\ref{iolargest})}                       \\ 
  &\ge&   c_2(\e)^{-1} \Big( \sum_{j \in \s_{k'}} |a_j|^2 \Big)^{1/2}
         - c_1(\e) \frac{r^{k_0 - k'}}{r-1} a
         \ \ \ \ \mbox{by (\ref{besselian})}                       \\
  &\ge&   \Big( c_2(\e)^{-1} - c_1(\e) (r - 1)^{-1} \Big) r^{k_0 - k'} a        
         \ \ \ \ \mbox{by (\ref{iolargest})}                        \\ 
  &\ge&   \Big( c_2(\e)^{-1} - c_1(\e) (r - 1)^{-1} \Big) a.
\end{eqnarray*}
If $r$ was chosen so that 
$c_2(\e)^{-1} - c_1(\e) (r - 1)^{-1}  >  c_2(\e)^{-1} / 2$, we are done.
The proof is complete.
\endproof

\noindent {\bf Remark 1. \ }
$C$ tends to $1$ as $\e \ra 1$. 
This is a consequence of a restriction theorem \cite{K-Tz} which we use 
in the following special case (see aslo \cite{B-Tz} Theorem 1.6).

\begin{theorem} (B.Kashin, L.Tzafriri).                                      \label{zeroes}
  Let $T$ be a linear operator in $l_2^n$ with $0$'s on the diagonal and $\|T\| = 1$. 
  Let $1/n  \le  \d  <  1$. 
  Then there exists a set $\s  \subset  \{1, \ldots, n\}$ with $|\s|  \ge  \d n / 4$
  for which 
  $$
  \| R_\s T R_\s \|  \le  c_5 \d^{1/2}.
  $$
\end{theorem}

First, Theorem \ref{quantitative} gives us a set of indices $\s_1$ with $|\s_1| \ge n/2$
such that the system $(x_j / \|x_j\|)_{j \in \s_1}$ is $c_6 c$-equivalent to the 
canonical vector basis of $l_2^{\s_1}$.
Let $\d = 1 - \e$ and $z_j  =  x_j / \|x_j\|$ for $j \in \s_1$.
Consider the linear operator $T$ in $l_2^{\s_1}$ which sends $e_j$ to $z_j$ for $j \in \s_1$.
Then the operator $T^* T - I$ has $0$'s on the diagonal and is of norm at most $2 c_6^2 c^2$. 
Applying Theorem \ref{zeroes} we get a set $\s \subset \s_1$ with $\s  \ge  \d |\s_1| / 4$
such that the following holds. 
For any sequence of scalars $(a_j)$
$$
\Big\| \Big\ll
(T^* T - I) \sum_{j \in \s} a_j e_j , \sum_{j \in \s} a_j e_j
\Big\rr \Big\|
\le  (2 c_6^2 c^2) c_5 \d^{1/2} = c_7 c^2 \d^{1/2}.
$$
Thus 
$$
\Big| \Big\ll  
\sum_{j \in \s} a_j z_j , \sum_{j \in \s} a_j z_j
\Big\rr
- \sum_{j \in \s} |a_j|^2 \Big|
\le  c_7 c^2 \d^{1/2}.
$$
Therefore the sequence $(z_j)_{j \in \s}$ is $g(\d)$-equivalent to $(e_j)_{j \in \s}$
for a function $g(\d)$ which tends to $1$ as $\d \ra 0$.
This proves Remark 1.

\vspace{0.5cm}
\noindent {\bf Remark 2. \ }
$h(\e)$ tends to infinity as $\e \ra \ 0$.
This is verified for the following tight frame $(x_j)_{1 \le j \le n+1}$, $n \ge 2$, 
considered by P.Casazza and O.Christensen in \cite{C-C2}:
\begin{eqnarray*}
  x_j     &=&  e_j - n^{-1} \sum_{j=1}^n e_j      \ \ \ \ \mbox{for $j = 1, \ldots, n;$}  \\
  x_{n+1} &=&  n^{-1/2} \sum_{j=1}^n e_j.
\end{eqnarray*}
Indeed, let $\s  \subset  \{1, \ldots, n\}$ be such that $|\s| > (1 - \e) n$
and the system $(x_j)_{j \in \s}$ is $M$-equivalent to an orthogonal basis.
By change of coordinates, the system $(x_j)_{1 \le j \le |\s|-1}$ must be $M$-equivalent
to an orthogonal basis as well. However,
$$
\Big\|  \sum_{j=1}^{|\s| - 1} x_j  \Big\|^2  \le  2 (\e n + 1)
$$
while $\|x_j\| \ge 1/2$ for all $j$.
Therefore $M$ can not be bounded independently of $n$ as $\e \ra 0$.
This proves Remark 2.

\section{Almost orthogonal subsequences of frames}              \label{basicsubs}

In this section we prove an infinite dimensional version of Theorem \ref{quantitative}.

\begin{theorem}                                                 \label{main}
  Given an $\e > 0$,
  every infinite dimensional frame has a subsequence 
  $(1 - \e)$-equivalent to an orthogonal basis of $l_2$. 
\end{theorem}

Given two sets $A$ and $B$ in $H$, we put by definition
$$
\theta (A, B)  =  \sup_{a \in A} \ \dist (a, B)
               =  \sup_{a \in A} \ \inf \{ \|a - b\| : b \in B \}.
$$

\begin{lemma}                                                 \label{distone}
  Let $(x_j)$ be a frame in an infinite-dimensional $H$.
  Let $A = \{ x_j / \|x_j\| \}$. 
  Then for any finite-dimensional subspace $E \subset H$
  $$
  \theta (A, E) = 1.
  $$
\end{lemma}

\proof
Let $z_j  =  x_j / \|x_j\|$ for all $j$.
Assume for the contrary that there is a  $\d < 1$ such that 
$$
\dist (z_j, E)  <  \d    \ \ \ \ \mbox{for all $j$}.
$$
Let $P$ be the orthogonal projection in $H$ onto $E$. Then 
$$
\|P z_j \|  >  \sqrt{1 - \d^2}   \ \ \ \ \mbox{for all $j$},
$$
so that 
\begin{equation}                                             \label{pxj}
  \|P x_j \|  \ge  \sqrt{1 - \d^2} \cdot \|x_j\|  \ \ \ \ \mbox{for all $j$}.
\end{equation}
Since $P$ is finite-dimensional, 
Lemma \ref{square} yields that the sequence $\|P x_j \|$ is square
summable. 
Then, by (\ref{pxj}), $\| x_j \|$ must be square summable, too. 
Thus $(x_j)$ is finite-dimensional. This contradiction completes the
proof.
\endproof

\begin{lemma}                                                 \label{stability}
  Let $\e_j$ be a sequence of quickly decreasing positive numbers
  ($2^{-j-1}$ will do).
  Let $(z_j)$ be a normalized sequence in $H$ such that 
  $$
  \ll z_i, z_j \rr  <  \e_j      \ \ \ \ \mbox{whenever $i < j$}.
  $$
  Then $(z_j)$ is equivalent to an orthonormal basis. 
\end{lemma}
The proof is simple.

\vspace{0.5cm}
{\bf Proof of Theorem  \ref{main}}. \ \
First note that, given an $\e > 0$, every subsequence equivalent to the
canonical vector basis of $l_2$ is weakly null, therefore has a subsequence 
which is $(1 - \e)$-equivalent to the canonical vector basis of $l_2$. 
Hence by Corollary \ref{frameistight} we may assume that our given 
frame $(x_j)$ is tight. 
Let $z_j  =  x_j / \|x_j\|$ for all $j$.
We will find a subsequence $(z_{j_k})$ equivalent to an orthogonal
basis by induction.
Put $j_1 = 1$. Let $j_1, \dots, j_{k-1}$ be defined 
and let $E = \span( z_{j_1}, \dots, z_{j_{k-1}} )$.
Choose $j_k$ from Lemma \ref{distone} so that
$$
\dist ( z_{j_k}, E )  >  1 - 2^{-2k}.
$$
Then it is easy to check that the constructed subsequence $(z_{j_k})$ 
satisfies the assumpiton of  Lemma \ref{stability}. This finishes the
proof.  
\endproof

\section{A frame not containing bases with brackets}                  
            \label{secwithoutbrackets}

\begin{definition}                                                    
            \label{basiswithbrackets}
  A sequence $(x_n)_{n=1}^\infty$ in a Banach space $X$ 
  is called a {\em basis with brackets} if there are numbers 
  $1  <  n_1  <  n_2  <  \ldots$ such that every vector $x \in X$
  admits a unique representation of the form
  $$
  x  =  \lim_j \sum_{n=1}^{n_j} a_n x_n, \ \ \ \ a_n \in \R.
  $$
\end{definition}

Clearly, every basis is a basis with brackets. 
The difference between bases and bases with brackets is that 
the latter require the convergence only of {\em some} partial sums 
in the representation.

The following lemma is known \cite{L-T}.

\begin{lemma}                                                         
            \label{triangle}
  Let $(x_n)_{n=1}^\infty$ be a basis with brackets, 
  and numbers $1  <  n_1  <  n_2  <  \ldots$ be as in 
  Definition \ref{basiswithbrackets}.
  Consider the projection $P_j$ onto $[x_n : n \le n_j]$
  parallel to $[x_n : n > n_j]$.
  Then $\sup_j \|P_j\| < \infty$.
\end{lemma}
Clearly, the converse also holds: if $\sup_j \|P_j\| < \infty$, 
for some sequence $1  <  n_1  <  n_2  <  \ldots$,
then $(x_n)$ is a basis with brackets.

\vspace{0.5cm}
In this section we prove

\begin{theorem}                                               \label{withoutbrackets}
  There exists a frame not containing bases with brackets.
\end{theorem}

\noindent Moreover, this frame is tight and have norms bounded from below.

\begin{lemma}                                                 \label{eitheror}
  There is an orthonormal basis $(z_j)$ in $l_2^n$ such that, 
  given any set $J \subset \{1, \ldots, n\}$, $|J| \ge n-2$, one has
  \begin{eqnarray*}
    \dist(e_1, [z_j : j \in J, j \ge j_0])  \le  4 / \sqrt{n}       
        \ \ \ \ & & \mbox{for $1   \le j_0  <  n/2$},               \\
    \dist(e_n, [z_j : j \in J, j  <  j_0])  \le  4 / \sqrt{n} 
        \ \ \ \ & & \mbox{for $n/2 \le j_0 \le n$}.                 
  \end{eqnarray*}   
\end{lemma}

\proof
By rotation, it is enough to find normalized vectors $v_1, v_2$ in
$l_2^n$
such that $\ll v_1, v_2 \rr  =  0$ 
and, given a set $J$ as in the hypothesis, 
\begin{eqnarray*}
  \dist(v_1, [e_j : j \in J, j \ge j_0])  \le  4 / \sqrt{n} 
      \ \ \ \ & & \mbox{for $1   \le j_0  <  n/2$},               \\
  \dist(v_2, [e_j : j \in J, j  <  j_0])  \le  4 / \sqrt{n} 
      \ \ \ \ & & \mbox{for $n/2 \le j_0 \le n$}.                 
\end{eqnarray*}      
Clearly, one may take
$$
v_1 = \halfn^{-1/2} \cdot 
  ( \underbrace{1, \ldots, 1}_\halfn, 0, \ldots, 0)   
\ \ \mbox{and} \ \
v_2 = \halfn^{-1/2} \cdot 
  ( 0, \ldots, 0, \underbrace{1, \ldots, 1}_\halfn).
$$
This completes the proof. 
\endproof

We will construct our frame $(x_j)$ by blocks $(x_j : j \in J(n))$,
where
$$
J(1) = \{1\},               \ \ 
J(2) = \{2, 3\},            \ \
J(3) = \{4, 5, 6\},         \ \
J(4) = \{7, 8, 9, 10\},     \ldots
$$
The supports of $x_j$'s from block $J(n)$ will lie in an interval
$I(n)$, where
$$
I(1) = \{1\},               \ \ 
I(2) = \{1, 2\},            \ \
I(3) = \{2, 3, 4\},         \ \
I(4) = \{4, 5, 6, 7\},     \ldots
$$
Let $i(n)$ be the first element in $I(n)$.

$$
\begin{array}{ccccccccccccc}
* & * & * &   &   &   &   &\lo&   &   &          \\
  & * & * & * & * & * &   &   &   &   &          \\
  &   &   & * & * & * &   &   &   &   &          \\
  &   &   & * & * & * & * & * & * & * &          \\  
  &   &   &   &   &   & * & * & * & * &          \\    
  &   &\lo&   &   &   & * & * & * & * &          \\    
  &   &   &   &   &   & * & * & * & * & \cdots  
 \end{array}
$$

The columns of this infinite matrix form the frame elements $x_j$, 
the asterisks marking their support.
Consider the shift operator $T_n : l_2^n \ra l_2$ 
which sends $(e_i)_{i=1}^n$ to $(e_i : i \in I(n))$.
Choose an orthonormal basis $(z_j : j \in J(n))$ in $l_2^n$ 
satisfying the conclusion of Lemma \ref{eitheror},
and define
$$
x_j  =  T_n z_j    \ \ \ \ \mbox{for $j \in J(n)$}.
$$

\begin{lemma}
  $(x_j)$ is a frame.
\end{lemma}  

\proof
Indeed, look at the rows in the picture, that is the vectors
$y_i  =  (x_1(i), x_2(i), \ldots)$.
Since the vectors $x_j$, $j \in J(n)$ are orthonormal for a fixed $n$, 
the vectors $y_i$ are orthogonal. Moreover, their norm is either equal 
to $2$ (if $i = i(n)$ for some $n$) or to $1$ (otherwise).
Now we pass again from the rows $y_i$ to the columns $x_j$. 
Lemma \ref{columnrow} yields that $(x_j)$ is a frame.
\endproof

Let $J$ be a set of positive integers such that 
the sequence $(x_j)_{j \in J}$ is complete in $l_2$. 
We shall prove that it is not a basis with brackets.

\begin{lemma}
  $|J(n) \cap J|  \ge  n - 2$ for every $n$.
\end{lemma}

\proof
Let $P$ be the orthogonal projection onto those $n-2$ coordinates in
$I(n)$
which don't belong to the other blocks $I(n_1)$, i.e. onto 
$[e_i : i \in I(n) \setminus \{i(n), i(n+1)\} ]$.
Thus $P$ sends to zero all $x_j$ with $j \not\in J(n)$. 
Hence $\im(P)  =  P( [x_j : j \in J(n) \cap J] )$.
Since $\im(P)$ is an $(n-2)$-dimensional space, the lemma follows.
\endproof

In the sequel we consider large blocks $J(n)$, i.e. with $n \ra \infty$.
Given a vector $v$ and a subspace $L$ in $l_2$ (both possibly
dependent on $n$),
we say that {\em $v$ is close to $L$} if $\dist(x, L) \le c / \sqrt{n}$.
Here $c$ is some absoulte constant, whose value may be different in
different occurences.

\begin{lemma}                                                         
                \label{onetwothree}
  1) $e_{i(n)}$   is close to $[x_j : j \in J(n-1) \cap J]$.

  2) $e_{i(n+1)}$ is close to $[x_j : j \in J(n+1) \cap J]$.
  
  3) Given a $j_0 \in J(n)$, 
     either $e_{i(n)}$   is close to $[x_j : j \in J(n) \cap J, j \ge
j_0]$,
     or     $e_{i(n+1)}$ is close to $[x_j : j \in J(n) \cap J, j  < 
j_0]$.
\end{lemma}

\proof
Note that $T_n$ sends $e_1$ to $e_{i(n)}$ and $e_n$ to $e_{i(n+1)}$.
Then all three statements of the lemma follow from Lemma \ref{eitheror}.
\endproof

The next and the last lemma, in tandem with Lemma \ref{triangle},
 completes the proof of Theorem \ref{withoutbrackets}.

\begin{lemma}
  For every $j_0 \in J(n)$ there is a normalized vector $x$ in $l_2$
  which is close to both subspaces
  $E = [x_j : j \in J, j \ge j_0]$ and 
  $F = [x_j : j \in J, j  <  j_0]$.
\end{lemma}

\proof
We make use of Lemma \ref{onetwothree}. 
By 3), we take either 
$x = e_{i(n)}$ to have $x$ close to $E$, or
$x = e_{i(n+1)}$ to have $x$ close to $F$.
In the first case $x$ is also close to $F$ by 2),
and in the second case $x$ is close to $E$ by 1). 
The proof is complete.
\endproof

\vspace{0.5cm}
A part of this work was accomplished when the author was visiting 
Friedrich-Schiller-Universit\"{a}t Jena. 
The author is grateful to M.Rudelson and P.Wojtaszczyk for helpful discussions,
and to V.Kadets for his constant encouragement.


{\small
\begin{thebibliography}{S 99}

\bibitem [A] {A}    A. Aldroubi, 
  {\em Portraits of frames},
  Proc. of the AMS 123 (1995), 1661--1668

\bibitem [B-Tz] {B-Tz}    J. Bourgain,  L. Tzafriri,
  {\em Invertibility of "large" submatrices with applications 
       to the geometry of Banach spaces and harmonic analysis},
  Israel J. Math. 57 (1987), 137--224

\bibitem [C1] {C1}    P. G. Casazza, 
  {\em Characterizing Hilbert space frames with the subframe property},
  Illinois J. Math. 41 (1997), 648--666

\bibitem [C2] {C2}    P. G. Casazza, 
  {\em Local theory of frames and Schauder bases for Hilbert space},
  Illinois J. Math., to appear (1999) 

\bibitem [C-C1] {C-C1}    P. G. Casazza, O. Christensen,
  {\em Hilbert space frames containing a Riesz basis 
       and Banach spaces which have no subspace isomorphic to $c\sb 0$},
  J. Math. Anal. Appl. 202 (1996), 940--950  

\bibitem [C-C2] {C-C2}    P. G. Casazza, O. Christensen,
  {\em Frames containing a Riesz basis 
       and preservation of this property under perturbations},
  SIAM J. Math. Anal. 29 (1998), 266--278     

\bibitem [D-S] {D-S}    R. J. Duffin, A. C. Schaeffer,
  {\em A class of nonharmonic Fourier series},
  Trans. of the AMS 72 (1952), 341--366

\bibitem [Ho] {Ho}    J. R. Holub,
  {\em Pre-frame operators, Besselian frames, and near-Riesz bases
       in Hilbert spaces},
  Proc. of the AMS 122 (1994), 779--785

\bibitem [He] {He}    C. Heil,
  {\em Wavelets and frames. Signal processing, Part I}, 
  147--160, IMA Vol. Math. Appl., 22, Springer, 1990

\bibitem [K-Tz] {K-Tz}    B. Kashin, L. Tzafriri, 
  {\em Some remarks on the restrictions of operators to coordinate subspaces},
  Preprint

\bibitem [L-T] {L-T}    J. Lindenstrauss, L. Tzafriri,
  {\em Classical Banach spaces},
  Springer, 1977

\bibitem [L] {L}    A. Lunin, 
  {\em On operator norms of submatrices},
  Matem. Zametki 45 (1989), 94--100

\bibitem [S] {S}    K. Seip, 
  {\em On the connection between exponential loss 
       and certain related  sequences in $L^2(-\pi, \pi)$},
  J. Funct. Anal. 130 (1995), 131--160     


\end{thebibliography}
}
\end{document}